\documentclass[11pt,twoside,english]{amsart}

\advance\oddsidemargin by -1.0cm
\advance\evensidemargin by -1.0cm
\advance\topmargin by -1.0cm
\usepackage{amssymb}
\usepackage{babel}
\usepackage{array}
\usepackage{url}

\theoremstyle{plain}

\newtheorem{lem}{Lemma}[section]
\newtheorem{thm}{Theorem}[section]
\newtheorem*{thm*}{Theorem}
\newtheorem{redthm}{Reduction Theorem}[section]

\theoremstyle{definition}

\newtheorem{NB}{Remark}[section]

\newcommand{\bdm}{\begin{displaymath}}
\newcommand{\edm}{\end{displaymath}}
\newcommand{\be}{\begin{equation}}
\newcommand{\ee}{\end{equation}}
\numberwithin{equation}{section}
\newcommand{\ba}[1]{\begin{array}{#1}}
\newcommand{\ea}{\end{array}}
\newcommand{\bea}[1][]{\begin{eqnarray#1}}
\newcommand{\eea}[1][]{\end{eqnarray#1}}
\newcommand{\btab}{\begin{tabular}}
\newcommand{\etab}{\end{tabular}}

\newcommand{\lan}{\left\langle}
\newcommand{\ran}{\right\rangle}

\newcommand{\ad}{\ensuremath{\mathrm{ad}}}

\newcommand{\C}{\ensuremath{\mathbb{C}}}
\newcommand{\R}{\ensuremath{\mathbb{R}}}

\newcommand{\eps}{\ensuremath{\varepsilon}}
\newcommand{\vphi}{\ensuremath{\varphi}}

\newcommand{\Scal}{\ensuremath{\mathrm{Scal}}}

\renewcommand{\ad}{\ensuremath{\mathrm{ad}\,}}

\newcommand{\grad}{\ensuremath{\mathrm{grad}\,}}

\newcommand{\n}{\ensuremath{\mathfrak{n}}}
\newcommand{\s}{\ensuremath{\mathfrak{s}}}

\newcommand{\hodge}{{*}}
\newcommand{\lie}[1]{\mathfrak{#1}}
\newcommand{\Lie}[1]{\textsl{#1}}
\newcommand{\G}{G\sb 2}

\DeclareMathOperator{\SL}{\Lie{SL}}
\DeclareMathOperator{\Spin}{\Lie{Spin}}
\DeclareMathOperator{\Sp}{\Lie{Sp}}
\DeclareMathOperator{\SU}{\Lie{SU}}

\newcommand{\g}{\lie{g}}

\newcommand{\pspm}{\psi^\pm}

\newcommand{\iso}{\cong}

\newcommand{\hook}{\lrcorner\,}

\newcommand{\sym}{\mathcal{S}}

\newcommand{\oo}{\omega^2}

\renewcommand{\leq}{\leqslant}

\newcommand{\q}{\quad}

\newcommand{\f}{\varphi}\newcommand{\ff}{\hodge\varphi}

\newcommand{\w}{\wedge}

\newcommand{\na}{\nabla}

\def\sideremark#1{\ifvmode\leavevmode\fi\vadjust{\vbox to0pt{\vss
 \hbox to 0pt{\hskip\hsize\hskip1em%
 \vbox{\hsize2cm\tiny\raggedright\pretolerance10000%
 \noindent #1\hfill}\hss}\vbox to8pt{\vfil}\vss}}}%

\begin{document}

\def\haken{\mathbin{\hbox to 6pt{%
                  \vrule height0.4pt width5pt depth0pt
                  \kern-.4pt
                  \vrule height6pt width0.4pt depth0pt\hss}}}
     \let \hook\intprod
\setcounter{equation}{0}
%
%
\thispagestyle{empty}
%
\title[Solvmanifolds with integrable and non-integrable
  $\G$ structures]
{Solvmanifolds with integrable and non-integrable $\boldsymbol{\G}$ structures}
%
%
%
\author{Ilka Agricola}
\author{Simon G. Chiossi}
\author{Anna Fino}
\address{\hspace{-5mm}
{\normalfont\ttfamily agricola@mathematik.hu-berlin.de}\newline
{\normalfont\ttfamily sgc@mathematik.hu-berlin.de}\newline
Institut f\"ur Mathematik \newline
Humboldt-Universit\"at zu Berlin\newline
Unter den Linden 6\newline
Sitz: John-von-Neumann-Haus, Adlershof\newline
D-10099 Berlin, Germany}
\address{\hspace{-5mm}
{\normalfont\ttfamily fino@dm.unito.it}\newline
Dipartimento di matematica\newline
Universit\`a di Torino \newline
Via Carlo Alberto 10\newline
10123 Torino, Italy}
%
\thanks{Supported by the SFB 647 "Space--Time--Matter" of the DFG, 
the Junior Research Group "Special Geometries in Mathematical Physics" of
the Volkswagen Foundation, GNSAGA of INdAM and MIUR}
\subjclass[2000]{Primary 53 C 25; Secondary 81 T 30}
\begin{abstract}
We show that a $7$-dimensional non-compact Ricci-flat Riemannian manifold
 with Riemannian holonomy $G_2$ can admit non-integrable $G_2$ structures
 of type $\R\oplus \sym^2_0(\R^7)\oplus \R^7$
in the sense of Fern\'andez and Gray.
 This relies on the construction of some $G_2$ solvmanifolds, whose
Levi-Civita connection is known to give a parallel spinor, admitting a 
$2$-parameter family of metric connections with non-zero 
skew-symmetric torsion that has
parallel spinors as well. The family turns out to be a deformation
of the Levi-Civita connection. This is in contrast with
 the case of compact scalar-flat Riemannian spin manifolds,
 where any metric connection with closed torsion admitting
 parallel spinors has to be torsion-free.

\end{abstract}
\maketitle
\pagestyle{headings}
%
%
%
\section{Introduction}\noindent
The study and explicit construction of Riemannian 
metrics with holonomy $G_2$ on non-compact manifolds of dimension
seven (called metrics with \emph{parallel} or \emph{integrable} $G_2$
  structure) has been an exciting area of differential
geometry since the pioneering work of Bryant and Salamon in the second
half of the eighties (cf.~\cite{Br}, \cite{BrS} and \cite{redbook}).
Mathematical elegance aside, 
these metrics have turned out to be an
important tool in superstring theory, since they are exact solutions of the 
common sector of type II string equations with vanishing $B$ field.

Independently of this development, the past years have
shown that non-integrable geometric structures such as almost hermitian
manifolds, contact structures or non-integrable $G_2$ and $\Spin(7)$ structures 
can be treated successfully with the powerful machinery of metric connections with
skew-symmetric torsion (see for example \cite{Friedrich&I1}, \cite{Agri},
\cite{AgFr1} and 
the literature cited therein). In physical applications, this torsion
is identified with a non-vanishing $B$ field (\cite{Strominger}, 
\cite{Gauntlett} and many more). The interaction between these
research lines was up to now limited to ``cone-type arguments'', i.\,e.~a
non-integrable structure on some manifold was used to
construct an integrable structure on a higher dimensional manifold
(like its cone, an so on). 
A natural question  is thus whether the same
Riemannian manifold $(M,g)$ can carry structures of both type 
\emph{simultaneously}. This appears to be a remarkable property. 
For example the projective space $\mathbb{CP}^3$
  with the well-known
 K\"ahler-Einstein structure and the nearly
K\"ahler one inherited from
triality does not satisfy this requirement. The metric underlying the nearly-K\"ahler 
structure is not the Fubini--Study one in fact, cf.~\cite{Eells-S:twistorial-harmonic}
and also \cite{BFGK}.

A spinor which is parallel with respect to a metric connection $\nabla$
(with or without torsion) forces its holonomy
to be a subgroup of the stabiliser of an algebraic spinor, and 
this is precisely $G_2$ in dimension $7$. In this particular dimension furthermore,
the converse statement also holds. The problem can therefore be 
reformulated as follows: \smallbreak

\textbf{Question.}
{\it Are there $7$-dimensional Riemannian 
manifolds with a parallel spinor for the Levi-Civita connection (rendering them
 Ricci-flat, in particular) also admitting a covariantly constant spinor 
for some other metric connection with skew-symmetric torsion? 
If yes, can the torsion connection be 
deformed into the Levi-Civita connection in such a way of preserving the parallel
spinor?  
} \smallbreak

\noindent
From the high energy physics' point of view a parallel spinor is interpreted as 
a supersymmetry transformation. Hence the physical problem behind the
 Question (which in fact motivated our investigations) is really whether a free 
``vacuum solution'' can also carry a non-vacuum supersymmetry, and how
the two are related. 

The case of a compact Riemannian manifold was treated in \cite{AgFr1}.
There, as a main application of the ``rescaled Schr\"odinger-Lichnerowicz 
formula'', one showed a rigidity theorem for compact manifolds
of non-positive scalar curvature. More precisely, 
\begin{thm*}
Suppose $(M^n, g, T)$ is
a compact, Riemannian spin
manifold of non-positive scalar curvature, $\Scal^g \leq 0$, and 
the $4$-form $dT$ acts on spinors as a non-positive endomorphism. Then
if there
exists a solution $\psi \neq 0$  of the equation
\bdm
\nabla^{T}_X\psi \ := \ \nabla^g_X \psi \, + \, (X \haken T) \cdot 
\psi \ = \ 0 \ ,
\edm
then the $3$-form and the scalar curvature vanish, $T = 0 = \Scal^g$,
and $\psi$ is parallel with respect to the Levi-Civita connection.
\end{thm*}

\noindent
This applies, in particular, to Calabi-Yau and
Joyce manifolds. These
are compact, Riemannian Ricci-flat manifolds of dimension $n=6,7,8$ with
(at least) one LC-parallel spinor field; under mild assumptions on the 
derivative of the torsion form $T$, they do not admit parallel spinors 
for any metric connection with $T\neq 0$. Since these manifolds have not been 
realized in any geometrically explicit way so far, harmonic or closed forms
are the natural candidates to be torsion forms on them.

The present paper deals with the non-compact case. Gibbons  
et al.~produced non-complete metrics
with holonomy $G_2$ in \cite{GLPS}.
Those metrics have the interesting feature,  among others, of admitting a $2$-step 
nilpotent isometry group $N$ acting on orbits of codimension one.
By \cite{Ch-F} such metrics are locally conformal to homogeneous metrics on
rank-one solvable extension of $N$, and the induced $SU(3)$ structure
on $N$ is half-flat. In the same paper all half-flat $SU(3)$
structures on $6$-dimensional nilpotent Lie groups whose
rank-one solvable extension is endowed with a conformally parallel 
$G_2$ structure were classified.
There are exactly six instances, which we considered in relation to
the problem posed.
It turns out that one of these manifolds provides a positive
answer to both questions
(Theorem \ref{main-result}), hence becoming the most interesting.
The wealth of parallel
spinors this manifold admits is organised into a continuous family
parametrised by the real projective line, plus a bunch of `isolated'
instances. To achieve this we proved a sort of `reduction' result
that allows to assume the spinors have an extremely simple block
form (Theorem \ref{thm-reduction}).
The Lie algebra associated to this
solvmanifold  has non-vanishing Lie brackets
\bdm
\ba{l}
 [e_i, e_7] = -\tfrac 35 m e_i,\ i = 1, 2, 5,\\[2pt]
 [e_j, e_7] = -\tfrac 65 m e_j,\ j = 3, 4, 6,\\[2pt]
 [e_1, e_5] = \tfrac 25 m e_3,\ [e_2, e_5] = \tfrac 25 m e_4,\ 
 [e_1, e_2] = \tfrac 25 m e_6.
\ea
\edm
The homogeneous metric it bears can be also seen as a $\G$ metric on
the product $\R\times \mathbb{T}$, where $ \mathbb{T}$ is the total
space of a $T^3$-bundle over another $3$-torus. 

Four metrics of the six only carry integrable 
$G_2$ structures (Theorem \ref{no-structure}), thus reproducing the
pattern of the compact situation, whilst the remaining one (example (4))
is singled out by complex 
solutions, a proper interpretation for which is still lacking (Theorem 
\ref{complex-structure}). Nevertheless, all the $G\sb 2$-metrics generated 
by these examples have a
physical relevance \cite{GLPS,Luest-T:AdS4-IIA}.
%
\section{General set-up}\noindent
%
The starting point of the present analysis is the classification of conformally parallel
$G_2$-manifolds on solvable Lie groups of \cite{Ch-F}, 
whose results we briefly summarise. We shall adopt a similar notation, 
except that the $3$-forms $\pspm$ have
become $\eta^\pm$, the conformal constant $m$ has changed sign to
$-m$, merely for aesthetic reasons, and the extension coefficients 
are now denoted by capital C's.
We shall also not distinguish between vectors and covectors.

\subsection{Round-up on $\boldsymbol{\G}$ solvable extensions} 
Consider a six-dimensional nilpotent 
Lie group $N$ with Lie algebra $\n$ endowed with an
invariant $SU(3)$ structure $(\omega,\eta^+)$, i.e. non-degenerate
$2$- and $3$-forms with stabilisers $\Sp(6,\R)$ and $\SL(3,\C)$
respectively. These define a
Riemannian metric with orthonormal basis $e_1,\ldots,e_6$ and an
orthogonal almost complex structure $J$.
Recall that $\ad_U(V)=[U,V]$ gives the
adjoint representation of a Lie algebra $\g$.
Pick the rank-one metric
solvable extension  $\s := \n \oplus \R e_7$, with $e_7 \perp\n$
a unit element, defined by $\ad_{e_7}$ as
non-singular self-adjoint derivation. 
The Lie bracket and 
inner product on $\s$ are, when restricted to $e_7^{\ \perp}$, 
precisely those of $\n$.

One is actually entitled to assume that there exists a unitary
basis $(e_1, \ldots, e_6)$ on $\n$ consisting of eigenvectors
of the derivation $\ad_{e_7}$ with
non-zero real eigenvalues $C_1, \ldots, C_6$.
In addition, all eigenvalues are positive integers without
common divisor, up to a rescaling of $e_7$
\cite{Heber:noncompact-Einstein, Will}.
Relatively to this basis of $\n$, the hermitian geometry of $N$
is prescribed  by
\bdm
\omega\ =\ e_{14}-e_{23}+e_{56},\quad
\eta^+ +i\eta^-= (e_1+ie_4)\wedge (e_2-ie_3)\wedge (e_5+i e_6),
\edm
The (non-integrable) $G_2$ structure inducing $g$
\bdm
\vphi\ :=\ \omega\wedge e_7+\eta^+ \ =\ e_{147}-e_{237}+e_{567} +
e_{125}+e_{136}+e_{246}-e_{345}
\edm
on the solvable Lie group $S$ corresponding to $\s$
is conformally parallel if and only if $\n$ is isomorphic
to one of the following:
\begin{enumerate}
\item  \q$ (0,0,e_{15}, 0,0,0)$, 
\item  \q$ (0,0,e_{15},e_{25},0,e_{12})$,
\item  \q$ (0,0,e_{15}-e_{46},0,0,0)$,
\item  \q$ (0,e_{45},-e_{15}-e_{46},0,0,0)$,
\item  \q$ (0,e_{45},e_{46},0,0,0)$,
\item  \q$ (0,e_{16}+e_{45},e_{15}-e_{46},0,0,0)$.
\end{enumerate}
The notation for Lie algebras is the usual differential one:
in (2) for instance, $e_{15}$ means $e_1\wedge e_5$ and the only
non-vanishing Lie brackets on $\n$ are
$[e_1,e_5]=-e_3,\ [e_2,e_5]=-e_4,\ [e_1,e_2]=-e_6$.
Throughout this article, the numeration shall respect the previous list.
\smallbreak

So the central issue here is the interplay of:
\begin{enumerate}
\item[(i)] the 6-dimensional manifold $(N, \omega, \eta^+)$;
\item[(ii)] the geometry of $S$ associated to the metric $g$ conformal to a 
  parallel one $\tilde g$;
\item[(iii)] the Ricci-flat metric $\tilde g$ on $S$ obtained by conformal change.
\end{enumerate}
We are mainly interested in the last structure, that is to say in the
incomplete metric $\tilde{g}$ with Riemannian holonomy 
contained in $G_2$.
We will show that in certain cases $\tilde g$ is induced by another  
$G_2$ structure, whose kind we describe. This helps to explain how
this non-integrable reduction is  related
to an integrable $G_2$ structure.\smallbreak

As a matter of fact, this is the expression for the
integrable $G_2$ structure on $(S,\tilde{g})$ with respect
to its (new) orthonormal basis as well. It is  known that  $\vphi$
defines a $\nabla^{\tilde{g}}$-parallel spinor $\Psi$ by
\begin{equation}
\label{G2-spinor}
\vphi(X,Y,Z)\ =\ \tfrac{1}{4}\langle X\cdot Y \cdot Z\cdot \Psi, \Psi\rangle,
\end{equation}
where dots denote Clifford multiplication and $\lan\,,\ran$ is the
scalar product in the spinor bundle. The constant $1/4$ is  arbitrary.

In terms of the seven-dimensional spin representation $\Delta_7$
used in \cite{AgFr1} (explicitly given in Section \ref{reduction}), the spinor $\Psi$
of \eqref{G2-spinor} has components
\be\label{eq:Psi}
\Psi\ =\ (0,0,0,0,1,1,-1,1).
\ee
Since $\Delta_7$ is the complexification
of a real representation, we assume all spinors to be real, unless
stated otherwise.
\subsection{Classification of $\boldsymbol\G$ structures}
The various $\G$-properties of $7$-manifolds $(S,\f)$ can be studied
using the approach of Fern\'andez and Gray \cite{Fernandez-G:G2},
i.e.~describing algebraically the four irreducible
$\G$-representations 
${\mathcal T}_i$ of the intrinsic torsion space  
\be\label{eq:G2torsion}
T^*S\otimes \g_2^\perp=\bigoplus_{i=1}^4{\mathcal T}_i\iso
\R\oplus\g_2\oplus\sym^2_0\R^7\oplus\R^7.
\ee
The first summand is merely spanned by $\f$, the second
denotes
the adjoint representation of $G_2$, whilst the third is 
the space of symmetric tensors on $\R^7$ with no trace.
The corresponding components $(\tau_1, \tau_2, \tau_3, \tau_4)$ of the intrinsic
torsion are uniquely defined differential forms such that 
\be\label{d-delta phi}
d\vphi\ =\ \tau_1 \ff+3\,\tau_4\wedge \vphi + 
*\tau_3,\qquad
\delta\vphi = -4 \hodge(\tau_4\wedge *\vphi)+
*(\tau_2\wedge \vphi),
\ee
with $\delta=-\hodge d \hodge$ the codifferential of forms, see \cite{Bryant:G2}.
For instance $\tau_1$ and the Lee form $\tau_4$ are given by
\bdm
\tau_1= g( d\f,\ff )/7 \q \text{and}\q \tau_4=- \hodge(\hodge d\f\wedge\f)/12.
\edm
Is is moreover known that $\tau_2 =0$ is equivalent to the existence 
of an affine connection $\tilde \nabla$ with skew-symmetric
torsion such that $\tilde \nabla \varphi = 0$
 \cite{Friedrich&I1}.
\smallbreak

What we mean by the ubiquitous and often abused terms 
{\sl integrable} (or {\sl parallel}) 
and {\sl non-integrable} is
\begin{itemize}
\item[] $\f$ is an integrable $\G$ structure $\iff$ $\tau_i=0$ for 
  $i=1,2,3,4$.\smallbreak

\item[] $\f$ is non-integrable $\iff$ one of the $\tau_i$'s at least 
survives, in which case the type of $\f$ is described by the non-zero
summands in \eqref{eq:G2torsion}.
\end{itemize}
This terminology is consistent with the landscape of general geometric
structures described in \cite{Fri2}.
\smallbreak

\noindent For example, a cosymplectic $\G$ structure $\f$ is characterised by the
equation $d\ff=0$, so it is non-integrable and has type 
${\mathcal T}_1\oplus{\mathcal T}_3\iso\R\oplus\sym^2_0\R^7$. 
The $\G$ structure of the previous page instead has type ${\mathcal
  T}_4$, as all $\tau_i$'s are zero except $\tau_4=me_7$.
\subsection{The Levi-Civita connection}
Let us sketch how one computes the torsion-free connection.
Denoting by $\hat{d}$ and $d$ the exterior differentials on $N$ and $S$,
the Maurer--Cartan equations for $\s=\n+\R e_7$ have the form
\bdm
d e_j\ =\ \hat{d}e_j+C_j e_{j7}\ \text{ for }j=1,\ldots,6\ \text{ and }
d e_7\ =\ 0.
\edm

\begin{table}
\begin{tabular}{|c|c|c|}\hline
Example & \hbox{$\n$ isomorphic to} & 
\hbox{eigenvalue type of $\ad_{e_7}$}\\[1mm]
\hline\hline
(1) & $ (0,0,e_{15}, 0,0,0)^{\vphantom{l}^{\vphantom{l}}}$ &
               $(2m/3, m, 4m/3, m,2m/3,m )$\\[1mm] \hline
(2) & $ (0,0,e_{15},e_{25},0,e_{12})^{\vphantom{l}^{\vphantom{l}}}$ &
               $(3m/5,3m/5,6m/5,6m/5,3m/5,6m/5)$\\[1mm]  \hline
(3) &  $ (0,0,e_{15}-e_{46},0,0,0)^{\vphantom{l}^{\vphantom{l}}}$ &
               $(3m/4, m,3m/2,3m/4,3m/4,3m/4)$\\[1mm]  \hline
(4) &  $ (0,e_{45},-e_{15}-e_{46},0,0,0)^{\vphantom{l}^{\vphantom{l}}}$ &
              $(4m/5,6m/5,7m/5,3m/5,3m/5,4m/5)$\\[1mm]  \hline
(5) & $ (0,e_{45},e_{46},0,0,0)^{\vphantom{l}^{\vphantom{l}}}$ &
               $(m,5m/4,5m/4,m/2,3m/4,3m/4)$\\[1mm]  \hline
(6) &  $ 
(0,e_{16}+e_{45},e_{15}-e_{46},0,0,0)^{\vphantom{l}^{\vphantom{l}}}$ &
               $(2m/3,4m/3,4m/3,2m/3,2m/3,2m/3)$\\[1mm]  \hline
\end{tabular}
\bigskip

\caption{The eigenvalue types and the underlying nilpotent Lie algebras $\n$.}
\end{table}
\noindent
The constant $m$ is real and positive, and it is important to remark
that each example is distinguished by a unique set of eigenvalues, as shown in Table 1.
A routine application of the Koszul formula yields the expression of
$\nabla^g$ on $S$ with respect to its
orthonormal basis $(e_1, \ldots, e_7)$. For instance it is not hard to see that
\bdm
\nabla^g_{e_i}e_7\ =\ C_ie_i,\quad \nabla^g_{e_i}e_i\ = \ -C_ie_i,\ 
\forall i\neq 7,\quad \nabla^g_{e_7}e_7\ =\ 0.
\edm
The new metric $\tilde{g}=e^{2f}g$ is determined by $df=me_7$.
The modified Levi-Civita connection 
can be computed through
\bdm
\nabla^{\tilde{g}}_X Y\ =\  \nabla^g_X Y + df(X)Y+df(Y)X-g(X,Y)\grad f,
\edm
so in particular
\be\label{eq:nabla-tildeg}
\ba{ll}
\nabla^{\tilde{g}}_{e_i}e_7=(C_i-1)e_i,& 
\nabla^{\tilde{g}}_{e_7}e_i=me_i,\\[3pt]
\nabla^{\tilde{g}}_{e_i}e_i=(1-C_i)e_i,& \nabla^{\tilde{g}}_{e_7}e_7=me_7
\ea
\ee
for all $i\neq 7$. 
The expression 
for the covariant derivatives
of the orthonormal basis $\tilde{e}_i:=e^{-f}e_i$ of 
$\tilde{g}$ can eventually be lifted to the spinor bundle.
We shall write $e_i$ instead of $\tilde{e}_i$ when no confusion arises. 
Therefore 
\begin{lem}\label{lemma:e7}
The derivatives of all vectors on $\n$ in the seventh direction are zero
\begin{center}
$\nabla^{\tilde{g}}_{e_7}U=0\q \text{for all}\ U\in\n$.
\end{center}
\end{lem}
\begin{proof}
This follows at once by conformally changing the relations in the 
second column of \eqref{eq:nabla-tildeg}.
\end{proof}
\noindent This will come handy in the next Section. 
%
\section{Reduction theorem for potential solutions}
\label{reduction}\noindent
Now we investigate whether the solvable Lie group 
$(S, \tilde g)$ admits a parallel
spinor for another metric connection with skew-symmetric torsion 
$T=\sum c_{\alpha\beta\gamma}\,e_{\alpha\beta\gamma}$. 
Instead of taking the most general $3$-form in dimension seven which
has 35 summands, we will make the Ansatz that $T$ be 
a linear combination
of the simple forms appearing in $\eta^+,\eta^-$ and $\omega\wedge e_7$.
Let $\Lambda^3_{11}(S)$ denote the subspace of $\Lambda^3(S)$ they
  span. 
Throughout the treatise
we shall take the spin representation $\Delta_7$ used in \cite{AgFr1, BFGK}:
\setlength{\extrarowheight}{4pt}
\bdm
\begin{array}{ll}
e_1 \, = \, +E_{18}+E_{27}-E_{36}-E_{45}, & e_2 \, =\, -E_{17}+E_{28}+E_{35}
-E_{46}, \\ 
e_3 \, =\, -E_{16}+E_{25}-E_{38}+E_{47}, & e_4\, =\, -E_{15}-E_{26}-E_{37}
-E_{48}, \\ 
e_5 \, =\,  -E_{13}-E_{24}+E_{57}+E_{68}, & e_6\, =\, +E_{14}-E_{23}-E_{58}
+E_{67}, \\ 
 e_7\, =\, +E_{12}-E_{34}-E_{56}+E_{78}, & \\ 
\end{array}
\edm
where $E_{ij}$ stands for the endomorphism of  $\R^7$ sending $e_i$ to 
$e_j$, $e_j$ to $-e_i$ and everything else to zero. 
Assuming then that the torsion looks like this
\bea[*]
T& =& 
c_{125}\,e_{125} + c_{136}\,e_{136} + c_{246}\,e_{246} +
c_{345}\,e_{345}+ \\ 
&&  c_{126}\,e_{126} + c_{346}\,e_{346}+ c_{135}\,e_{135} +  c_{245}\,e_{245} + \\
&& c_{147}\,e_{147} + c_{567}\,e_{567} + c_{237}\,e_{237},
\eea[*]
and denoting by $\haken$ the interior product, one infers that 
Clifford multiplication by $e_i\haken T$
 has -- as an endomorphism -- the block structure  
$\bigl(\begin{smallmatrix} 0 & * \\   * & 0 \end{smallmatrix}\bigr)$ 
for any $i$. This is
particularly interesting when $i=7$ in the light of Lemma \ref{lemma:e7}.
It allows one to determine the structure of elements in $\ker
(\nabla^{\tilde{g}}_{e_7} + e_7\haken T)= \ker ( e_7\haken T)$ without
too much effort.
Clearly only the coefficients
$c_{147}, c_{237}$ and $c_{567}$ of $T$ are involved. 
\begin{redthm}\label{thm-reduction}
%
For $T\in\Lambda^3_{11}(S)$ a non-trivial element annihilated by $e_7\haken T$
is a linear combination of {\rm upper block} reduced forms
\begin{itemize}
\item[(A)] $\psi=(a,b,c,d,0,0,0,0)$ with
$c_{567}=0$, $c_{147}= -c_{237}$,\smallbreak
\item[(B)] $\psi=(a,b, -\eps\,a, \eps\,b,0,0,0,0)$ 
and $c_{147}= -c_{237}+\eps c_{567}$ with $\eps=\pm 1$,\smallbreak
\end{itemize}
or {\rm lower block} reduced forms
\begin{itemize}

\item[(C)] $\psi=(0,0,0,0,e,f,g,h)$ with
$c_{567}=0$, $c_{147}= +c_{237}$,\smallbreak
\item[(D)] $\psi=(0,0,0,0,e,f,\eps\,e,-\eps\,f)$ 
and $c_{147}= +c_{237}+\eps c_{567}$ with $\eps=\pm 1$.
\end{itemize}
\end{redthm}
\begin{NB}
Notice that the cases are not mutually exclusive: for example
if $c_{567}=0$, (B) is a special case of (A) as (D) is of (C).
\end{NB}
\noindent
In conclusion, one can always assume that a spinor has such a block structure,
with the coefficients $c_{147}, c_{237}, c_{567}$ subjected to
one addtional linear constraint.
%
\section{Families of real solutions}
\noindent
%
The solvable extension of Example (2) 
is equipped with a Ricci-flat 
metric with Riemannian holonomy equal to $G_2$, implying
that there exists a unique $\nabla^{\tilde g}$-parallel spinor
$\psi$. In terms of the
endomorphisms $E_{ij}$, the Levi-Civita
connection on the tangent bundle has components
\bdm
\ba{ll}
\nabla^{\tilde{g}}_{e_1} = -\tfrac 15 me^{-f} (2E_{17}+E_{35}-E_{26}),&
\nabla^{\tilde{g}}_{e_2} = -\tfrac 15 me^{-f} (E_{16}+2E_{27}+ E_{45}),\\[3pt]
\nabla^{\tilde{g}}_{e_3} = -\tfrac 15 me^{-f} (E_{15}-E_{37}),&
\nabla^{\tilde{g}}_{e_4} = -\tfrac 15 me^{-f} (E_{25}-E_{47}),\\[3pt]
\nabla^{\tilde{g}}_{e_5} = -\tfrac 15 me^{-f} (E_{13}+E_{24}+2E_{57}),&
\nabla^{\tilde{g}}_{e_6} = -\tfrac 15 me^{-f} (E_{12}-E_{67}),
\ea
\edm
and $\nabla^{\tilde{g}}_{e_7}=0$. 

Let us now study the existence of solutions $\psi\neq 0$ of the equation
$\nabla_{e_i}^T \psi= 0$, where by definition
\be\label{parallelflach2}
\nabla_{e_i}^T\psi = \nabla_{e_i}^{\tilde g}\psi+ (e_i \haken T)\cdot\psi.
\ee
\begin{thm}\label{main-result}
The  equation $\nabla^T\psi=0$
admits precisely $7$ solutions for some
$T\in \Lambda^3_{11}(S)$, namely:
\begin{itemize}
\item[a)] A two-parameter family of pairs $(T_{r,s},\psi_{r,s})$ such that
$\nabla^{T_{r,s}}\psi_{r,s}=0$;\smallbreak
\item[b)] Six `isolated' solutions occuring in pairs, 
$(T^\eps_i,\psi^\eps_i)$ for $i=1,2,3$ and $\eps=\pm$.
\end{itemize}
All $G_2$ structures admit one parallel spinor, and for 
\begin{itemize}
\item[] $ |r|\neq |s|$:\q $\f_{r,s}$ is of general type 
  $\R\oplus \sym^2_0\R^7\oplus \R^7$,\smallbreak
\item[] $r=s$:\q $\vphi_{r,r}$ is parallel, the torsion $T_{r,r}=0$ and 
  $\psi_{r,r}$ is a multiple of $\psi$.\smallbreak
\item[] $r=-s$:\q the $\G$ type of $\f_{-s,s}$ has no $\R$-term.
\end{itemize}
\end{thm}
\begin{proof}
By Reduction Theorem \ref{thm-reduction} we can treat cases
(A)-(D) separately. This yields the following possibilities.\medbreak

\noindent   
\textbf{Solution a)}. Set
\be\label{eq:l-m}
\lambda_{r,s}\ =\ \frac{r^2-s^2}{2(r^2+s^2)},\quad \mu_{r,s}\ =\ 
\frac{(r-s)^2}{r^2+s^2}
\ee
and
\bdm
\psi_{r,s}\ =\ (0,0,0,0,r,s,-r,s).
\edm
The spinor $\psi_{r,s}$ is parallel with respect to
the connection $\nabla^{r,s}:=\nabla^{T_{r,s}}$
determined by 
\bea[*]
T_{r,s} & =&  
-\tfrac{1}{10}me^{-f}\left[\lambda_{r,s}(\eta^+ - 6\,e_{125})+
\mu_{r,s}(\eta^- +3\,e_{346}) \right].
\eea[*]
Notice that this family of $3$-forms contains no terms in $e_7$. Furthermore,
$\lambda_{r,s}=\lambda_{cr,cs}$ and $\mu_{r,s}=\mu_{cr,cs}$ for any
real constant $c\neq 0$, reflecting
the fact that any multiple of $\psi_{r,s}$ is again parallel for the 
connection with the same torsion form. The $G_2$
structure corresponding to $\psi_{r,s}$ is
\be\label{eq:f_rs}
\vphi_{r,s}\ =\ rs\,\eta^+ + \frac{1}{2}(s^2-r^2)\,\eta^- + 
\frac{1}{2}(s^2+r^2)\,\omega\wedge e_7.
\ee
It is by now clear why taking $r=\pm s$ plays a special role, for $T_{r,s}$ 
and $\vphi_{r,s}$ both simplify. 
The type of $\vphi_{r,s}$ is determined once one computes its 
differential and codifferential.
Recall that from the covariant derivative of a $3$-form $\xi=e_{ijk}$,
\bdm
\nabla_X(e_{ijk})\ =\ (\nabla_X e_i)\wedge e_{jk}+ e_i\wedge
(\nabla_x e_j)\wedge e_k+ e_{ij}\wedge (\nabla_X e_k),
\edm
one obtains $d$ and $\delta$ by
\bdm
d \xi(X_0,\ldots,X_3) = \sum_{i=0}^3(-1)^i(\nabla_{X_i}\xi)(X_0,\ldots,
\hat{X}_i,\ldots, X_3),\qquad
\delta\xi\ =\ -\sum_i e_i\haken \nabla_{e_i}\xi.
\edm
The result of these lengthy calcultations is given in Table $2$.

For $r=s$ all components $\tau_i$ of the intrinsic torsion vanish, 
since $\vphi_{r,s}$ is integrable. 
By construction $\tau_4$ is proportional to $e_7$, with constant $c$
resulting from the discussion.
In general, \eqref{d-delta phi} gives
\bdm
\ba{l}
d\f_{r,s}\ =\ \tfrac{s^2-r^2}2 d\eta^-+\tfrac{(s-r)^2}2 d\eta^+\\[3pt]
\delta\vphi_{r,s}\ =\ -\frac{1}{5} me^{-f}(r-s)^2\omega
\q\text{ and }\q -4 \hodge(c\,e_7\wedge *\vphi_{r,s})\ =\
-2c(r^2+s^2)\omega.
\ea
\edm
This implies that $c=\tfrac 1{10} me^{-f}\mu_{r,s} \neq 0$  
 for $r\neq s$ and $\tau_2$ is identically zero, as one 
expects. 
As for $\tau_1=-\tfrac 3{10} me^{-f}(r^2-s^2)(2r^2+2s^2-rs)$, one 
sees it also vanishes for $r=-s$, since $e_{1257}$ does not appear 
in $*\vphi_{r,s}$.
The $4$-form 
\bea[*]
\hodge\tau_3 & = & 
 -\tfrac 35 m(s-r)^2e_{1257}+\tfrac
3{10}m(s^2-r^2)(\eta^-+2e_{346})\w e_7+ \tfrac 1{10}m(s^2-r^2)\oo+\\
&& 
\tfrac 3{10}m(s^2-r^2)(2s^2+2r^2-sr)\bigl(-rs\eta^-\w e_7+\tfrac{s^2-r^2}2
\,\eta^+\w e_7 -\tfrac{s^2+r^2}4\,\oo\bigr) -\\
&& \tfrac
3{10}m\,\tfrac{(r-s)^2}{s^2+r^2}\bigl(rs\,\eta^+\w e_7+\tfrac{(s^2-r^2)}2
(\eta^+-2e_{126})\w e_7\bigr)
\eea[*]
is never zero for $r\neq s$, instead.
\begin{table}
\setlength{\extrarowheight}{4pt}
\begin{tabular}{|c|c|c|c|}
\hline
\hbox{Form} & \hbox{differential $d$}  & \hbox{Hodge $*$}  &
\hbox{codifferential $\delta $} \\ 
\hline\hline
$e_{125}$ & $-\frac{6}{5}me^{-f} e_{1257}$ & $e_{3467}$ &
$\frac{2}{5}me^{-f}\omega$ \\[1mm] \hline
$e_{136}$ & $ 0$ & $e_{2457} $  & $0$ \\[1mm] \hline
$e_{246}$ & $ 0$ & $e_{1357} $  & $0$ \\[1mm] \hline
$e_{345}$ & $ 0$ & $e_{1267} $  & $0$ \\[1mm] \hline\hline
$e_{126}$ & $-\frac{3}{5}me^{-f} e_{1267} $ & $-e_{3457} $  & $0$ \\[1mm] \hline
$e_{135}$ & $-\frac{3}{5}me^{-f} e_{1357} $ & $-e_{2467} $  & $0$ \\[1mm] \hline
$e_{245}$ & $-\frac{3}{5}me^{-f} e_{2457} $ & $-e_{1367} $  & $0$ \\[1mm] \hline
$e_{346}$ & $ \frac{1}{5}me^{-f}
    \oo+\frac{3}{5}me^{-f}e_{3467} $ & $-e_{1257} $  &
$0$ \\[1mm]\hline\hline
$e_{147}$ & $\frac{2}{5}me^{-f}e_{1257} $ & $e_{2356} $  &
$-\frac{2}{5}me^{-f}e_{14} $ \\[1mm] \hline
$e_{237}$ & $\frac{2}{5}me^{-f}e_{1257} $ & $e_{1456} $  &
$-\frac{2}{5}me^{-f}e_{23} $ \\[1mm] \hline
$e_{567}$ & $\frac{2}{5}me^{-f}e_{1257} $ & $e_{1234} $  &
$-\frac{2}{5}me^{-f}e_{56} $ \\[1mm] \hline
\end{tabular}

\bigskip

\caption{Derivatives of the simple forms spanning
  $\Lambda^3_{11}(S)$.}
\end{table}

\bigskip\noindent
\textbf{Solution b)}. The isolated solutions occur in pairs labelled $\pm$,
 basically corresponding to the choice of sign for $\eps$ in 
the Reduction Theorem.
The first couple consists of the spinors 
\bdm
\psi^+_1= (0,1,0,-1,0,0,0,0) \q\text{and}\q \psi^-_1=(1,0,1,0,0,0,0,0)
\edm
(denoted $\psi^\eps_1$ with $\eps=\pm$) and the $3$-forms
\bea[*]
T^\eps_1 &=&
-\frac{me^{-f}}{10} \left[ \frac{\eps}{2}(\eta^+ +4\,e_{125}-2e_{246})+
\frac{1}{3}(\eta^- -2e_{135}- e_{346}) -\frac{2\eps}{3}
(\omega-e_{23})\wedge e_7 \right].
\eea[*]
The additional relation on the $c_{ij7}$'s reads
$c_{147}=-c_{567}-c_{237}$. 
Via equation \eqref{G2-spinor} 
the characteristic form is 
\bdm
2\, \vphi^\eps_1\ =\ \eps(e_{126}+e_{135}-e_{245}+e_{346}) -
e_{147}-e_{567}- e_{237}. 
\edm
\medbreak
\noindent The second pair of solutions 
gives spinors  
\bdm
\psi^+_2= (0,1,0,1,0,0,0,0) \q\text{and}\q
\psi^-_2=(1,0,-1,0,0,0,0,0),
\edm
 together with the torsion
\bea[*]
T^\eps_2 &=& 
-\frac{me^{-f}}{10} \left[ \frac{\eps}{2}(\eta^+ +4\,e_{125}-2e_{136})+
\frac{1}{3}(\eta^- + 2e_{245}- e_{346}) -\frac{2\eps}{3}
(\omega+e_{14})\wedge e_7  \right].
\eea[*]
The underlying relation is $c_{147}=c_{567}-c_{237}$. The
characteristic $3$-form is
\bdm
2\, \vphi^\eps_2\ =\ \eps(-e_{126}+e_{135}-e_{245}-e_{346}) +
e_{147}-e_{567}+ e_{237}.
\edm
\medbreak
\noindent For the last pair, the spinors are lower block 
\bdm
\psi^+_3=(0,0,0,0,1,0,1,0) \q\text{and}\q
\psi^-_3=(0,0,0,0,0,1,0,-1).
\edm
The torsion $3$-form is then
\bea[*]
T^\eps_3 &=& 
-\frac{me^{-f}}{10} \left[ \frac{1}{2}(\eta^+ +4\,e_{125}+2e_{345})+
\frac{1}{3}(\eta^- -2e_{126}- e_{346}) -\frac{2\eps}{3}
(\omega+e_{56})\wedge e_7  \right].
\eea[*]
In this case the equation $c_{147}=c_{567}+c_{237}$ holds. Now 
the characteristic $3$-form is 
\bdm
2\, \vphi^\eps_3\ =\ \eps(e_{126}+e_{135}+e_{245}-e_{346}) -
e_{147}+e_{567}+ e_{237}.
\edm
In all cases it is not hard to check that $\vphi^\eps_i$ have type
$\R\oplus \sym^2_0(\R^7)\oplus \R^7$.  \qedhere
\end{proof}
\begin{NB}
The family of $\G$ structures \eqref{eq:f_rs} depends upon the 
two homogeneous parameters \eqref{eq:l-m}, 
or if one prefers on the projective coordinate $w=r/s$. 
In fact $\lambda=\lambda_{r,s}, \mu=\mu_{r,s}$ lie on the 
ellipse $(\mu-1)^2+4\lambda^2-1=0$
in the $(\lambda,\mu)$-plane.
The extremal points $w=\infty,0$ correspond to
$\f_{r,0}=\tfrac {r^2}2(-\eta^-+\omega\w e_7)$ and 
$\f_{0,s}=\tfrac {s^2}2(+\eta^-+\omega\w e_7)$, where $\eta^+$ is
missing. 
Similarly, the origin of $\R^2$ is $\f_{r,r}=r^2(\eta^++\omega\w e_7)$ 
whilst $w=-1$ produces the
form $\f_{r,-r}=r^2(-\eta^++\omega\w e_7)$, and the roles of
$\eta^\pm$ are swapped. 
It is interesting perhaps to notice that each
$w$ on the conic $\mathbb{RP}^1$ corresponds to a specific choice
of $3$-form in the canonical bundle of $N$, 
and does not touch significantly the term $\omega\w e_7$.
\end{NB}
%
\section{The other examples}\noindent
%
The solvmanifolds extending numbers (1), (3), (5), and (6) admit no 
non-trivial solutions to \eqref{parallelflach2}, whereas (4) yields only 
\emph{complex} solutions.
We quickly gather the results, writing in particular the 
Levi-Civita connection.
%
\subsection*{Example (1)}
\noindent
%
In many respects, this example is the closest to the Riemannian flat case
$\R^7$. Although trivially Ricci-flat, Euclidean space admits no parallel
spinors for a connection with non-vanishing skew-symmetric torsion
\cite{AgFr1}. Here a similar result holds. The Riemannian holonomy
reduces to $\SU(2)\subset G_2$, and only three components of the
LC connection survive, precisely
\bdm
\nabla^{\tilde{g}}_{e_1} =\ -\tfrac 13 me^{-f} (E_{17}+E_{35}),\quad
\nabla^{\tilde{g}}_{e_3} =\ -\tfrac 13 me^{-f} (E_{15}-E_{37}),\quad
\nabla^{\tilde{g}}_{e_5} =\ -\tfrac 13 me^{-f} (E_{13}+E_{57}),
\edm
and the four $\nabla^{\tilde{g}}$-parallel spinors are
\bdm
(1,1,0,0,0,0,0,0),\quad (0,0,-1,1,0,0,0,0),\quad (0,0,0,0,1,1,0,0),\quad
(0,0,0,0,0,0,-1,1).
\edm
%
\subsection*{Example (3) }
\noindent
%
The Levi-Civita connection on the tangent bundle is given by
\bdm 
\ba{ll}
\na^{\tilde{g}}_{e_1} =\ -\tfrac{me^{-f}}{4} (E_{17}+E_{35}), &
\na^{\tilde{g}}_{e_2}\ =\ 0,\\[3pt]
\na^{\tilde{g}}_{e_3} =\ -\tfrac{me^{-f}}{4} (E_{15}-2E_{37}-
E_{46}), &
\na^{\tilde{g}}_{e_4}\ =\ -\tfrac{me^{-f}}{4} (-E_{36}+E_{47}),\\[3pt]
\na^{\tilde{g}}_{e_5} =\ -\tfrac{me^{-f}}{4} (E_{13}+E_{57}), &
\na^{\tilde{g}}_{e_6}\ =\ -\tfrac{me^{-f}}{4} (E_{34}+ E_{67}).
\ea\edm
It has holonomy group $SU(3)$, so $\Psi$  of \eqref{eq:Psi} pairs up with a second 
LC-parallel spinor 
 $(1,1,1,-1,0,0,0,0)$.
%
\subsection*{Example (5)}
\noindent
The Levi-Civita
connection is given by
\bdm 
\ba{ll}
\na^{\tilde{g}}_{e_1} =\ 0, &
\na^{\tilde{g}}_{e_2}\ =\ -\tfrac{me^{-f}}{4} (-E_{27}-E_{45}),\\[3pt]
\na^{\tilde{g}}_{e_3} =\ -\tfrac{me^{-f}}{4} (-E_{37}- E_{46}), &
\na^{\tilde{g}}_{e_4}\ =\ -\tfrac{me^{-f}}{4} (-E_{25}-E_{36}+2E_{47}),\\[3pt]
\na^{\tilde{g}}_{e_5} =\ -\tfrac{me^{-f}}{4} (E_{24}+E_{57}), &
\na^{\tilde{g}}_{e_6}\ =\ -\tfrac{me^{-f}}{4} (E_{34}+ E_{67}).
\ea\edm
The holonomy is $SU(3)$, hence there exists another 
$\nabla^{\tilde{g}}$-parallel spinor besides 
$\Psi$, namely
$(-1,1,1,1,0,0,0,0)$.

\subsection*{Example (6) }
\noindent
%
The Levi-Civita
connection on the tangent bundle is given by
\bdm
\ba{ll}
\nabla^{\tilde{g}}_{e_1} = -\tfrac{me^{-f}}{6} (2E_{17}+E_{35}-E_{26}),& \quad
\nabla^{\tilde{g}}_{e_2}\ =\ -\tfrac{me^{-f}}{6} (-E_{16}-2E_{27}- E_{45}),\\
\nabla^{\tilde{g}}_{e_3} = -\tfrac{me^{-f}}{6} (E_{15}-2E_{37}-E_{46}),& \quad
\nabla^{\tilde{g}}_{e_4}\ =\ -\tfrac{me^{-f}}{6} (-E_{25}-E_{36}+2E_{47}),\\
\nabla^{\tilde{g}}_{e_5} = -\tfrac{me^{-f}}{6} (E_{13}+E_{24}+2E_{57}),& \quad
\nabla^{\tilde{g}}_{e_6}\ =\ -\tfrac{me^{-f}}{6} (-E_{12}+E_{34}+ 2\,E_{67}),
\ea\edm
This manifold has full holonomy $G_2$. 
Then again
\begin{thm}\label{no-structure}
Let $(S,\tilde{g})$ be one of the solvmanifolds
$(1)$, $(3)$, $(5)$, or $(6)$.
If there exists a non-zero spinor $\psi$ solving $\nabla^T\psi=0$ for some
$T\in \Lambda^3_{11}(S)$, then $T=0$ and $\psi$ is a linear combination of
the given $\na^{\tilde{g}}$-parallel spinors.
\end{thm}
\begin{proof}
By the Reduction Theorem one can assume that  $\psi$
has a block structure. Considering cases (A)-(D) separately tells 
that there are no solutions except for $T=0$. 
\end{proof}
%
\section{Complex solutions}
\noindent 
The Riemannian connection of the manifold (4) reads
\bdm
\ba{ll}
\nabla^{\tilde{g}}_{e_1} = -\tfrac{me^{-f}}{5} (E_{17}+E_{35}),&\quad
\nabla^{\tilde{g}}_{e_2}\ =\ -\tfrac{me^{-f}}{5} (-E_{27}- E_{45}),\\
\nabla^{\tilde{g}}_{e_3} = -\tfrac{me^{-f}}{5} (E_{15}-2E_{37}-E_{46}),&\quad
\nabla^{\tilde{g}}_{e_4}\ =\ -\tfrac{me^{-f}}{5} (-E_{25}-E_{36}+2E_{47}),\\
\nabla^{\tilde{g}}_{e_5} = -\tfrac{me^{-f}}{5} (E_{13}+E_{24}+2E_{57}),&\quad
\nabla^{\tilde{g}}_{e_6}\ =\ -\tfrac{me^{-f}}{5} (E_{34}+E_{67}).
\ea
\edm
There are similarities with the 
Levi-Civita expression relative to example (2), although the two
solvmanifolds are not isometric. 

It is rather curious to be 
in presence of \emph{complex} solutions. Though one is usually
interested in real spinors and differential forms, complex coefficients  might as well
be relevant for other considerations. As in proof of Theorem \ref{main-result}, by 
the reduction process of \ref{thm-reduction} we can consider
the occurring cases one by one.

\begin{thm}\label{complex-structure}
Let $(S,\tilde{g})$ be the solvmanifold of example $(4)$. 
If there exists a non-zero spinor $\psi$ satisfying Equation 
\eqref{parallelflach2} for some
$T\in \Lambda^3_{11}(S)$ and all $i=1,\ldots,7$, then:
\begin{itemize}
\item[(a)] $\psi $ is a multiple of 
$\big(1+2i\eps\sqrt{2},\, 3,\, 1+2i\eps\sqrt{2} ,\, -3, 0,0,0,0\big)$
and 
\bea[*]
T & = & \tfrac 23\,[-2 e_{126}+ e_{135}-4 e_{245}+e_{346}] +
 i \eps\sqrt{2}\,[e_{125}+e_{136}+e_{246}+e_{345}]\\
&+& \tfrac 23\,i \eps\sqrt{2}\,[-e_{147}-e_{567}+ 2e_{237}],\q \text{or}
\eea[*]
\item[(b)] $\psi $ is a multiple of 
$\big(3, \,-1+2i\eps\sqrt{2},\, -3,\, -1+2i\eps\sqrt{2} , 0,0,0,0\big)$
and 
\bea[*]
T & = & \tfrac 23\,[e_{126}- e_{135}+4 e_{245}-2e_{346}] +
 i \eps\sqrt{2}\,[-e_{125}+e_{136}+e_{246}-e_{345}]\\
&+& \tfrac 23\,i \eps\sqrt{2}\,[-e_{147}+e_{567}+ 2e_{237}],\q \text{or}
\eea[*]
\item[(c)] $\psi $ is a multiple of 
$\big(0,0,0,0,\, 1+2i\eps\sqrt{2},\, 3,\, 1+2i\eps\sqrt{2} ,\, -3\big)$
and 
\bea[*]
T & = & \tfrac 23\,[e_{126}-2 e_{135}+4 e_{245}-e_{346}] +
 i \eps\sqrt{2}\,[e_{125}+e_{136}-e_{246}-e_{345}]\\
&+& \tfrac 23\,i \eps\sqrt{2}\,[e_{147}-e_{567}+ 2e_{237}].
\eea[*]
\end{itemize}
Above $\eps$ is $1$ or $-1$ and stems from the solution
of a quadratic equation.
\end{thm} 
\begin{NB}
In Strominger's model of superstring theory (\cite{Strominger},
\cite{Friedrich&I1}), the contraction $T(i,j):=\sum_{m,n}T_{imn}T_{jmn}$ 
appears as a relevant term, essentially the torsion contribution to the
Ricci tensor. A question of interest is then whether the term is real
for the complex solutions above. Now $T(i,j)$ is a real number,
possibly zero, apart when $e_i=\pm J(e_j)$
%
\begin{itemize}
\item[(a)] \ $T(1,4)\ =\ T(2,3)\ =\ T(5,6)\ = -8/3\sqrt{2}i\eps$\smallbreak
\item[(b)] \ $T(1,4)\ =\ -T(2,3)\ =\ T(5,6)\ = 8/3\sqrt{2}i\eps$\smallbreak
\item[(c)] \ $T(1,4)\ =\ -T(2,3)\ =\ -T(5,6)\ = 8/3\sqrt{2}i\eps$.
\end{itemize}
%
There seems to be no physical meaning for these solutions in the
models currently under investigation.
\end{NB}

\noindent 
It is tempting to pursue the same analysis without the
assumption that {\sl all} coefficients $C_j$ of the solvable extension $S$ 
be non-zero, which is important only in connection to the existence
of Einstein metrics on $S$ \cite{Heber:noncompact-Einstein}. 
With hindsight, we reasonably
expect to find metrics with holonomy strictly contained in $\G$, so
the developed technique might furnish many parallel spinors.
%
%

\end{document}